\documentclass[11pt,letterpaper]{amsart}
\pdfoutput=1

\usepackage{amssymb}
\usepackage{amsmath}
\usepackage{amscd}
\usepackage{longtable}
\usepackage{float}
\usepackage{bigstrut}

\usepackage{graphicx}
\usepackage{multirow}
\usepackage[all,cmtip]{xy}
\def\@begintheorem#1#2{\par\bgroup{\sc #1\ #2. }\it \ignorespaces}
\def\@opargbegintheorem#1#2#3{\par\bgroup{\sc #1\ #2\ (#3) . }\it \ignorespaces}
\def\@endtheorem{\egroup}

\newtheorem{theorem}{Theorem}[section]
\newtheorem{corollary}[theorem]{Corollary}
\newtheorem{lemma}[theorem]{Lemma}
\newtheorem{proposition}[theorem]{Proposition}
\newtheorem{example}[theorem]{Example}
\newtheorem{definition}[theorem]{Definition}
\newtheorem{remark}{Remark}[section]

\newcommand{\bt}[1]{\begin{theorem}\label{#1}}
\newcommand{\bc}[1]{\begin{corollary}\label{#1}}
\newcommand{\bl}[1]{\begin{lemma}\label{#1}}
\newcommand{\bp}[1]{\begin{proposition}\label{#1}}
\newcommand{\be}[1]{\begin{example}\label{#1}}
\newcommand{\bd}[1]{\begin{definition}\label{#1}}
\newcommand{\brem}[1]{\begin{remark}\label{#1}}
\newcommand{\bpr}{\par\noindent{\it Proof}: \ignorespaces}
\newcommand{\beq} {\begin{eqnarray}}

\newcommand{\et}{\end{theorem}}
\newcommand{\ec}{\end{corollary}}
\newcommand{\el}{\end{lemma}}
\newcommand{\ep}{\end{proposition}}
\newcommand{\ee}{\end{example}}
\newcommand{\ed}{\end{definition}}
\newcommand{\erem}{\end{remark}}
\newcommand{\eeq}{\end{eqnarray}}

\newcommand{\mepr}{{\ \ \ \ \ \ \vbox{\hrule\hbox{%
\vrule height1.3ex\hskip0.8ex\vrule}\hrule}}}

\def \I {{\mathcal I}}
\def \r {{\mathcal R}}
\def \q {{\mathcal Q}}

\title{\bf Orbit closures of representations of source-sink Dynkin quivers}
\author{\bf Kavita Sutar}
\date{}
\address{Department of Mathematics, Northeastern University}
\email{sutar.k@husky.neu.edu}

\begin{document}

\begin{abstract}
We use the geometric technique described in \cite{MR1988690, 2011arXiv1111.1179S} to calculate the resolution of orbit closures of representations of Dynkin quivers with every vertex being source or sink. We use this resolution to derive the normality of such orbit closures. As a consequence we obtain the normality of certain orbit closures of type $E$.
\end{abstract}
\maketitle

\section{Introduction} \label{sec1}

Fix an algebraically closed field $K$. A \emph{quiver} is a pair $Q=(Q_0,Q_1)$ where $Q_0$ is a set of vertices and $Q_1$ is a set of arrows. We denote by $\hat{Q}$ the underying graph of a quiver $Q$. 
 We use the notation $ta \stackrel{a}{\rightarrow} ha$ for arrows in $Q$. A \emph{source-sink quiver} will mean a quiver with every vertex being either a source or a sink. \par

A \emph{representation} $((V_i)_{i \in Q_0}, (V(a))_{a \in Q_1})$ of $Q$ is an assignment of finite dimensional $K$-vector spaces $V_i$ to every vertex $i \in Q_0$ and $K$-linear maps $V_{ta} \stackrel{V(a)}{\rightarrow} V_{ha}$ to every arrow $a \in Q_1$.  The \emph{dimension vector} of a representation $((V_x)_{x \in Q_0}, (V(a))_{a \in Q_1})$ is defined as the function $\underline{d}:Q_0 \longrightarrow \mathbb{Z}$ given by $\underline{d}(x)=$ dim $V_x$. The Euler form $E_Q$ of a quiver $Q$ is a quadratic map $E_Q:\mathbb{Z}^{|Q_0|} \rightarrow \mathbb{Z}$ given by 
\[(\alpha_i)_{i \in Q_0} \longmapsto \sum_{x \in Q_0} \alpha_x^2 -  \sum_{{a \in Q_1}|\stackrel{a}{ta \rightarrow ha}} \alpha_{ta} \alpha_{ha}  \]
\noindent Given two representations $V=((V_i)_{i \in Q_0}, (V(a))_{a \in Q_1})$ and $W=((W_i)_{i \in Q_0}, (W(a))_{a \in Q_1})$ of $Q$, a morphism $\Phi:V \rightarrow W$ is a collection of $K$-linear maps $\phi_i:V_i \rightarrow W_i$ such that for every $a \in Q_1$, the square
   \[
    \xymatrix{
    V_{ta} \ar[d]^{\phi_{ta}} \ar[r]^{V(a)} & V_{ha} \ar[d]^{\phi_{ha}}\\
    W_{ta} \ar[r]^{W(a)} &W_{ha}
    } 
    \]
commutes. \par

With this definition of morphisms, the collection of all representations of a quiver $Q$ (over $K$) forms a category which we denote by Rep$_K(Q)$. Given a quiver $Q$, one can define its path algebra $KQ$ as the $K$-algebra generated by the paths in $Q$. It is known that $KQ$ is an associative algebra and is finite dimensional if and only if $Q$ is finite and has no oriented cycles. An important result in the theory of representation theory of associative algebras asserts that for $Q$ being a finite, connected, acyclic quiver, there is an equivalence of categories Mod $KQ$ and Rep$_K(Q)$ (refer \cite{MR2197389} for details).\par

The \emph{representation space} $Rep(Q,\underline{d})$ of a quiver $Q$ is the collection of all representations of $Q$ of fixed dimension vector $\underline{d}$. We view $Rep(Q,\underline{d})$ as the set $\displaystyle{\prod_{a \in Q_1} \hbox{Hom}(K^{d_{ta}}, K^{d_{ha}})}$. Thus, $Rep(Q,\underline{d})$ is a finite dimensional $K$-vector space with an affine structure. The algebraic group $\prod_{x \in Q_0} GL_{d(x)}(K)$ acts on $Rep(Q,\underline{d})$ with action given by
\[((g_x)_{x \in Q_0}, V ) \longmapsto (g_{ha} V(a) g_{ta}^{-1})_{a \in Q_1}\] For $V \in Rep(Q,\underline{d})$, let $\overline{O}_V$ denote the closure of an orbit $O_V$. Then $\overline{O}_V$ is a subvariety of $Rep(Q,\underline{d})$. These varieties are in some sense a generalization of determinantal varieties.\par

In this paper we calculate a minimal free resolution of such varieties. The method is a generalization of Lascoux's calculation of the resolutions of determinantal varieties \cite{MR520233}. Using this resolution we draw conclusions about geometric properties like normality, rational singularities etc. of the orbit closures. The study of these geometric properties has been an interesting field of research during the past decade \cite{Z}. The question of normality of $\overline{O}_V$ in the case of Dynkin quivers of type $A$ and $D$ has been investigated by Abeasis, Del Fra and Kraft \cite{MR626958}, Bobinski-Zwara \cite{MR1967381,MR1885816} and Lakshmibai-Magyar \cite{MR1635873}. The question of normality for orbit closures corresponding to Dynkin quivers of type $E$ is open. We obtain a result (Corollary \ref{cor3.4}) which answers this question for a class of orbit closures corresponding to source-sink Dynkin quivers.

Section \ref{sec2} contains some preliminaries. We present the main results in Section \ref{sec3} and some examples of the calculation in Section \ref{sec4}.

\section{Preliminaries}\label{sec2}
 
\subsection{The geometric technique}
 We briefly sketch the geometric technique used for our calculations. For details we refer to \cite{MR1988690} (Chapter 5). \par
Let $X$ be affine space of dimension $N$ and $Y \subset X$ be a subvariety. Let $\mathcal{V}$ be a projective space of dimension $m$. Then $X \times \mathcal{V} \stackrel{p}{\rightarrow} \mathcal{V}$ is the trivial vector bundle of rank $N$, we denote it by $\mathcal{E}$. Suppose $Z \subset X \times  \mathcal{V}$ such that $Z$ is the total space of a subbundle $\mathcal{S}$ of $\mathcal{E}$.  Let $q: X \times \mathcal{V} \stackrel{q}{\rightarrow} X$ be the first projection. The Koszul complex resolving the structure sheaf $\mathcal{O}_Z$ of $Z$ can be used to calculate a complex $F_{\bullet}$. If $Z$ is a desingularization of $Y$ then the module $q_*(\mathcal{O}_Z)$ is the normalization of $K[Y]$. Under additional conditions, $F_{\bullet}$ gives a minimal free resolution of the coordinate ring of the normalization $\hat{Y}$ of $Y$. \par
Consider the exact sequence of vector bundles over $\mathcal{V}$
\[0 \longrightarrow \mathcal{S} \longrightarrow \mathcal{E} \longrightarrow \mathcal{T} \longrightarrow 0\]

Let $\xi=\mathcal{T}^*$. The Koszul complex 
\[\mathcal{K}(\xi)_{\bullet} : 0 \rightarrow \bigwedge^t(p^*\xi) \rightarrow \dots \rightarrow \bigwedge^2(p^*\xi) \rightarrow p^*\xi \rightarrow 0  \]
resolves the structure sheaf $\mathcal{O}_Z$ as $\mathcal{O}_{X \times \mathcal{V}}$-module. The complex  ${\mathbf F}_{\bullet}$ is obtained as a pushforward of the above Koszul complex and the terms of ${\mathbf F}_{\bullet}$ are calculated using Theorem \ref{thm2.1} and Bott's algorithm. 

\bt{thm2.1}(Basic theorem \cite{MR1988690}) The terms of the complex ${\mathbf F}_{\bullet}$ are given by
 \[ \textbf{F}_i = \bigoplus_{j \geq 0}H^j(\mathcal{V}, \bigwedge^{i+j}\xi)\otimes A[-i-j] \]
\et 

The following theorem tells us how we can use the complex ${\mathbf F}_{\bullet}$ to draw conclusions about the varieties in question.
\bt{thm2.3}\cite{MR1988690} With notation as above,
\begin{itemize}
\item[(1)] If ${\mathbf F}_i=0$ for $i < 0$ then ${\mathbf F}_{\bullet}$ is a finite free resolution of the normalization of $K[Y]$.
\item[(2)] If $\textbf{F}_i = 0$ for $i < 0$ and $\textbf{F}_0 = A$ then the variety $Y$ is normal and has rational singularities. 
\end{itemize}
\et

\subsection{Desingularization}    
To calculate the complex ${\mathbf F}_{\bullet}$ in our case we consider a desingularizaton $Z$ of an orbit closure $\overline{O}_V$ given by Reineke's construction \cite{MR1985731}. We describe this construction briefly here.

Let $Q$ be a Dynkin quiver and let $AR(Q)$ denote its corresponding Auslander-Reiten quiver. Let $\I$ be a partition of the Auslander-Reiten quiver satisfying:
  \begin{enumerate}
     \item $Ext^1_Q(X_{\alpha}, X_{\beta})=0$ for all $\alpha, \beta \in \I_t$ for $t=1,\cdots ,s$.
     \item $Hom_Q(X_{\beta}, X_{\alpha})=0= Ext^1_Q(X_{\alpha}, X_{\beta})$ for all $\alpha \in \I_t, \beta \in \I_u,t < u$
   \end{enumerate} 

Such a partition of  $AR(Q)$ exists because the category of finite-dimensional representations is directed; in particular, we can choose a sectional tilting module and let $\I_t$ be its Coxeter translates.  We fix a partition $\I_*$ of $AR(Q)$. Then the indecomposable representations $X_{\alpha}$ are the vertices of $AR(Q)$. For a representation $V= \oplus_{\alpha \in R^+}m_{\alpha} X_{\alpha}$, we define representations \[ V_{(t)}:=\oplus_{\alpha \in \I_t}m_{\alpha} X_{\alpha}, ~~~~ t=1,\cdots,s \] Then $V=V_{(1)}\oplus \cdots \oplus V_{(s)}$. Let $\underline{d}_t =$ dim $V_{(t)}$. We consider the incidence variety \[Z_{\I_*, V} \subset \prod_{x \in Q_0} Flag(d_s(x),d_{s-1}(x)+d_s(x), \cdots , d_2(x)+ \cdots +d_s(x), K^{d(x)} ) \times Rep_K(Q,\underline{d}) \] defined as
  \begin{equation}
  Z_{\I_*, V} = \{((R_s(x) \subset R_{s-1}(x) \subset \cdots \subset R_2(x) \subset K^{d(x)}), V) ~|~ \forall a \in Q_1, \forall t,   
  ~~V_a(R_t(ta)) \subset R_t(ha) \}
  \label{Z}
  \end{equation}
\bt{Thm 1}(Reineke \cite{MR1985731}) Let $Q$ be a Dynkin quiver, $\I_*$ a directed partition of $R^+$. Then the second projection \[ q:Z_{\I_*, V} \longrightarrow Rep_K(Q, \underline{d})\] makes $Z_{\I_*, V}$ a desingularization of the orbit closure $\overline{O}_V$. More precisely, the image of $q$ equals $\overline{O}_V$ and $q$ is a proper birational isomorphism of $Z_{\I_*, V}$ and $\overline{O}_V$.
\et   

In this case we say that $Z=Z_{\I_*, V}$ is a $(s-1)$-step desingularization. For our calculations, we restrict to orbit closures admitting a $1$-step desingularization. Then the vector bundle $\xi$ is 
     \begin{equation}
     \xi=\bigoplus_{a \in Q_1} \r_{ta}\otimes \q_{ha}^* 
     \label{xi} 
     \end{equation}
We calculate   
     \begin{eqnarray} \nonumber
     \bigwedge^t \xi &=& \bigoplus_{\sum_{a \in Q_1}{k_a}=t }\bigwedge^{k_a}(\r_{ta}\otimes \q_{ha}^* )\\ \nonumber
     &=& \bigoplus_{\sum_{a \in Q_1}{|\lambda(a)|=t}}\left[\bigotimes_{a \in Q_1}  
     S_{\lambda(a)}\r_{ta}\otimes S_{\lambda(a)'}\q_{ha}^* \right]~~~~~~~~~~ \hbox{(by Cauchy's formula)}\\
     &=&\bigoplus_{\sum_{a \in Q_1}{|\lambda(a)|=t}}\bigotimes_{x \in Q_0}\left[ \bigotimes_{a \in Q_1 | 
     ta=x}S_{\lambda(a)}\r_x \otimes \bigotimes_{a \in Q_1 | ha=x}S_{\lambda(a)'}\q_x^* \right] 
     \label{extpowxi}
     \end{eqnarray}      
  The term in (\ref{extpowxi}) is a vector bundle over $\mathcal{V}$ to which we associate a weight in the following manner. First some notation: if $x$ is a vertex with more than one incoming or outgoing vertices then the corresponding term in the right hand side of Equation (\ref{extpowxi}) is calculated using the Littlewood-Richardson rule. In such cases we will use the shorthand notation $\lambda(a_1a_2)$ to denote a Young tableau occuring in the Littlewood-Richardson product of Young tableaux $\lambda(a_1)$ and $\lambda(a_2)$. So for example, if the arrows $a_1, a_2, a_3$ are all the outgoing arrows from source $x$, then a summand of $S_{\lambda(a_1)}R_x \otimes S_{\lambda(a_2)}R_x \otimes S_{\lambda(a_3)}R_x$ will be denoted by $S_{\lambda(a_1a_2a_3)}R_x$. \par

Also we use notation $-\lambda$ for the non-increasing sequence consisting of terms $\lambda$ written with a minus sign and in reverse order (for example if $\lambda=(3~3~2~1)$, then $-\lambda=(-1~-2~-3~-3)$). With this notation we can describe the associated weight as follows: for $x \in Q_0$ let $a_1, a_2, \cdots a_k$ be all the outgoing arrows and $b_1, b_2, \cdots b_l$ be all the incoming arrows at $x$. Then the weight associated to the summand corresponding to vertex $x$ in Equation (\ref{extpowxi}) is
\[ (\underbrace{ -\lambda(b_1 b_2 \dots b_l)'}_{dim Q_x},~ \underbrace{\lambda(a_1 a_2 \dots a_k)}_{dim R_x}) \]

We apply Bott's algorithm to these weights to calculate the terms in Theorem \ref{thm2.1}.

We remark that we can calculate the complex ${\mathbf F}_{\bullet}$ starting with an incidence variety $Z(\beta \subset \alpha)$ introduced by Schofield in \cite{MR1162487}. These are defined as follows. Let $X=Rep(Q,\beta+\gamma)$ and $\mathcal{V}=\prod_{x \in Q_0} Gr(\beta_x,\beta_x+\gamma_x)$.  \[Z(Q,\beta \subset \beta+\gamma) \subset Rep(Q,\beta+\gamma) \times \prod_{x \in Q_0} Gr(\beta_x,\beta_x+\gamma_x)\] is defined as the collection of quiver representations of dimension vector $\alpha=\beta+\gamma$ together with a subrepresentation of dimension vector $\beta$. Thus
\[Z(Q,\beta \subset \beta+\gamma)=\{(V,R)\in Rep(Q,\beta+\gamma) \times \prod_{x \in Q_0} Gr(\beta_x,\beta_x+\gamma_x)~|~\forall a \in Q_1, V_a(R_{ta}) \subset R_{ha} \}\] 
 In the case of Dynkin quivers, the variety $Y=q(Z(Q,\beta \subset \beta+\gamma))$ is an orbit closure: $Z$ is irreducible implies $Y$ is irreducible and since there are only finitely many orbits in case of Dynkin quivers, we have that $Y$ must be an orbit closure.  In general however, it is not known whether $Y$ is an orbit closure. 
 
\section{Main results}\label{sec3}

First we have some results involving Young tableaux. \par
Notation: A partition $\lambda=(\lambda_1, \lambda_2, \cdots, \lambda_n)$ is a non-decreasing sequence of non-negative integers. The Young diagram corresponding to partition $\lambda$ consists of $\lambda_i$ boxes in the $i$th row. The conjugate partition $\lambda'$ is the partition $(\lambda_1', \lambda_2', \cdots, \lambda_m')$ where $\lambda_j'$ is the number of boxes in the $j$th column. We will denote the last row of a Young tableau $\lambda$ by $\lambda_{last}$. \par
The following lemma is an easy exercise in counting boxes- 
\bl{lemma3.1}
Let $\lambda$ be a Young tableau. Then for all $a$ and $b$,
 \[ \lambda_1+\lambda_2+ \cdots +\lambda_a \leq ab+(\lambda'_{b+1}+\cdots+\lambda'_{last}).  \]
\el
\bpr We consider three cases: \par
  \begin{itemize}
  \item[Case (1)] $\lambda'_{b+1}=a$. Then
        \[\lambda_1+\lambda_2+ \cdots \lambda_a= ab+\lambda'_{b+1}+\cdots+\lambda'_{last}\]
    
  \item[Case (2)] $\lambda'_{b+1} > a$. In this case $\lambda'_{b+1},\lambda'_{b+2},\cdots \lambda'_{last}$ contribute more boxes so that
   \[\lambda_1+\lambda_2+ \cdots +\lambda'_a\leq ab+\lambda'_{b+1}+\cdots +\lambda'_{last} \]
  \item[Case (3)] $\lambda'_{b+1}<a$. Here the rectangle $ab$ contributes more boxes, so that
  \[\lambda_1+\lambda_2+ \cdots +\lambda_a \leq ab+\lambda'_{b+1}+\cdots + \lambda'_{last}\]
  \end{itemize}
\begin{flushright} \mepr \end{flushright}

 By symmetry we also have for all $a$ and $b$:
     \begin{equation} \lambda'_1+\lambda'_2+ \cdots +\lambda'_a \leq ab+(\lambda_{b+1}+\cdots+\lambda_{last}) \label{ineq2}
     \end{equation}   

The next lemma is one of the well known Horn-type inequalities for triples of partitions \cite{MR1685640}. 
\bl{lemma3.2} Suppose $\nu$ is one of the partitions occuring in Littlewood-Richardson product of $\lambda$ and $\mu$. Then \[\nu_1+\nu_2+\cdots+\nu_k \leq (\lambda_1+\lambda_2+\cdots+\lambda_k)+(\mu_1+\mu_2+\cdots+\mu_k).\]
\el
\vskip0.7cm 
Let $Q=(Q_0, Q_1)$ be an source-sink Dynkin quiver. Fix a representation $V$ of $Q$. Let $(\lambda(a))_{a \in Q_1}$ be a $|Q_1|$-tuple of partitions. Consider the  variety $Z$ obtained as a $1$-step desingularization of $\overline{O}_V$ . When calculating the resolution ${\mathbf F}_{\bullet}$ of $q_*(\mathcal{O}_Z)$ we are concerned with the difference 
\[D(\underline{\lambda}):=\sum_{a \in Q_1}{|\lambda(a)|}-N\] 
where $N$ is the total number of Bott exchanges required in the process of obtaining a partition from the weights described below. Our main theorem is an inequality involving the above difference and $E_Q$. It is a generalization of the inequality obtained for $D(\lambda, \mu,\nu)$ in \cite[Proposition 4.4]{ 2011arXiv1111.1179S}. \par
Let $Q' \subset Q_0$ be the set of all source vertices and $Q'' \subset Q_0$ be the set of all sink vertices. Let $\lambda(a)$ be a non-increasing sequence associated to every arrow $a \in Q_1$. With this notation, the exterior power $\bigwedge^t \xi$ in Equation (\ref{extpowxi}) can be viewed as
  \begin{equation} \bigwedge^t \xi = \bigoplus_{\sum_{a \in Q_1}{|\lambda(a)|=t}}\left[(\bigotimes_{x \in Q'} \bigotimes_{a \in Q_1 | 
     ta=x}S_{\lambda(a)}\r_x )\otimes (\bigotimes_{x \in Q''}\bigotimes_{a \in Q_1 | ha=x}S_{\lambda(a)'}\q_x^*) \right]        
  \label{extpowxi1}
  \end{equation}
Thus we have one summand for every $|Q_1|$-tuple of non-increasing sequences $(\lambda(a))_{a \in Q_1}$. It will be useful to let this tuple of partitions also stand for the summand it corresponds to. \par
If $x$ is a vertex with more than one incoming or outgoing vertices then the corresponding term in the right hand side of Equation (\ref{extpowxi1}) is calculated using the Littlewood-Richardson rule for tensor products. Recall that we use the notation $\lambda(a_1a_2)$ to denote a Young tableau occuring in the Littlewood-Richardson product of Young tableaux $\lambda(a_1)$ and $\lambda(a_2)$. \par

To calculate the resolution ${\mathbf F}_{\bullet}$ we associate a weight to each summand of $\bigwedge^t \xi$. Each summand consists of tensor products of terms of the form $S_{\lambda(-)}R_x$ (for $x \in Q'$) and $S_{\lambda(-)}Q_x^*$ (for $x\in Q''$). If $x \in Q'$ the  associated sequence is $(0^{\gamma_x}, \lambda(a_1 a_2 \dots a_k))$ where $a_1, a_2, \cdots a_k$ are all the outgoing arrows at $x$; if $x \in Q''$, the sequence is $(-\lambda(b_1 b_2 \dots b_l)', 0^{\beta_x})$ where $b_1, b_2, \cdots b_l$ are all the incoming arrows at $x$. 
We can now state the main theorem.  
\bt{thm3.3} With notation as above,
\[D(\underline{\lambda}) \geq E_Q\]
\et
\bpr To calculate $D(\underline{\lambda})$ we apply Bott's algorithm to the weights described above and count the total number of exchanges $N$. There is one weight associated to every vertex; let $N_x$ denote the number of Bott exchanges at vertex $x$.  \par
If $x$ is a source, the weight at $x$ is of the form $(0^{\gamma_x}, \lambda(I_x))$ where $I_x=a_{i_1}a_{i_2}\dots a_{i_k}$ such that $a_{i_1}, a_{i_2},\dots a_{i_k}$ are all the arrows incident at $x$. Then $N_x=\gamma_x  u_x$ where $u_x$ is the largest number such that $\lambda(I)_{u_x} - \gamma_x \geq u_x$.\par
Similarly, if $y$ is a sink, then weight at $y$ is of the form $(-\lambda(J_y)', 0^{\beta_y})$, where $J_y=b_{j_1}b_{j_2}\dots b_{j_l}$ such that $b_{j_1}, b_{j_2},\dots b_{j_l}$ are all the arrows incident at $y$. In this case $N_y=\beta_y u_y$ where $u_y$ is the largest number such that $-\lambda(J)_{u_y} + \beta_y \leq u_y$. Thus\par
 \begin{equation} N=\sum_{x\in Q' } N_x+\sum_{y \in Q''}N_y = \sum_{x \in Q'} \gamma_x u_x + \sum_{y \in Q''} \beta_y u_y
 \label{N}
 \end{equation}
Note that if $u_x$ is the largest number such that $\lambda(I)_{u_x} - \gamma_x \geq u_x$ then \[\lambda(I_x)_1 \geq \lambda(I_x)_2 \geq \cdots \lambda(I_x)_{u_x} \geq \gamma_x + u_x \] implies
  \begin{equation}
  \lambda(I_x)_1+\lambda(I_x)_2+\cdots+ \lambda(I_x)_{u_x} \geq u_x(\gamma_x+u_x) = u_x^2+\gamma_xu_x 
  \label{term1}
  \end{equation}
For similar reasons we have
  \begin{equation}
  \lambda(J_y)'_1+\lambda(J_y)'_2+\cdots + \lambda(J_y)'_{u_y} \geq u_y(\beta_y+u_y) = u_y^2+\beta_yu_y 
  \label{term2}
  \end{equation} 
  
On the other hand we have by Lemma \ref{lemma3.2} that 
  \[\lambda(I_x)_1+\cdots \lambda(I_x)_{u_x} \leq \sum_{\stackrel{a_{i_k}}{x \longrightarrow}}(\lambda(a_{i_k})_1+\cdots + \lambda(a_{i_k})_{u_x})\]
  
Combining this with Inequality (\ref{term1}) gives 
 \begin{equation} 
\sum_{\stackrel{a_{i_k}}{x \longrightarrow}}(\lambda(a_{i_k})_1+\cdots + \lambda(a_{i_k})_{u_x}) \geq u_x^2+\gamma_xu_x 
 \label{term3} 
 \end{equation}
for every pair $(x, I_x)$ with $x \in Q'$. \par

Similarly  
\[\lambda(J_y)'_1+\cdots+ \lambda(J_y)'_{u_y} \leq \sum_{\stackrel{b_{j_k}}{ \longrightarrow y}}(\lambda(b_{j_k})_1+\cdots + \lambda(b_{j_k})_{u_y})\]
together with Inequality (\ref{term2}) implies
    \begin{equation}
\sum_{\stackrel{b_{j_k}}{ \longrightarrow y}}(\lambda(b_{j_k})_1+\cdots + \lambda(b_{j_k})_{u_y}) \geq u_y^2+\beta_yu_y
  \label{term4}
  \end{equation}
for every pair $(y, J_y)$ with $y \in Q''$.\par
  
Using Lemma \ref{lemma3.1} we get a further upper bound on the right hand side terms of Inequality (\ref{term4}): if $b_{j_k}$ is an arrow from $x_k$ to $y$ then   
\[u_{x_k} u_{y}+\lambda(b_{j_k})_{u_{x_k}+1}+\lambda(b_{j_k})_{u_{x_k}+2}+\cdots+\lambda(b_{j_k})_{last} \geq \lambda(b_{j_k})'_1+\cdots +\lambda(b_{j_k})'_{u_{y}}  \]
for every $k=1,2,\dots ,l$. So for every pair $(y,J_y)$ we get inequalities 
\begin{equation}
\sum_{\stackrel{b_{j_k}}{x_k \rightarrow y}} (u_{x_k}u_y +\lambda(b_{j_k})_{u_{x_k}+1}+\lambda(b_{j_k})_{u_{x_k}+2}+\cdots+\lambda(b_{j_k})_{last} )\geq u_y^2+\beta_yu_y
\label{term5}
\end{equation}

Adding the inequalities in (\ref{term3}) and (\ref{term5}) for all pairs $(x,I_x)_{x \in Q'}$ and $(y,I_y)_{y \in Q''}$, we get 
\begin{equation}
\sum_{a \in Q_1}|\lambda(a)| + \sum_{\stackrel{a}{x \rightarrow y}} u_x u_y \geq \sum_{x \in Q'}(u_x^2+\gamma_xu_x)+\sum_{y \in Q''}  (u_y^2+\beta_yu_y)=\sum_{x \in Q_0}u_x^2 - N
\end{equation}
which means
\begin{equation}
\sum_{a \in Q_1}|\lambda(a)|-N \geq \sum_{x \in Q_0}u_x^2 -  \sum_{\stackrel{a}{x \rightarrow y}} u_x u_y 
\end{equation}
 
\begin{flushright} \mepr \end{flushright}

\bc{cor3.4} Let $Q$ be a Dynkin quiver with source-sink orientation, $V$ be a representation of $Q$ such that the orbit closure $\overline{O}_V$ admits a $1$-step desingularization $Z$. Then $\overline{O}_V$ is normal and has rational singularities.
\ec
\bpr $Q$ is Dynkin implies $E_Q > 0$. Theorem \ref{thm3.3} implies that the terms ${\mathbf F}_i$ of the resolution ${\mathbf F}_{\bullet}$ are zero for $i<0$ and $F_0=A$. By Theorem \ref{thm2.3} it follows that the orbit closure is normal and has rational singularities.
\qed 

\bc{cor3.5} Let $Q$ be an extended Dynkin quiver with source-sink orientation. If $V$ is a representation of $Q$ such that the orbit closure $\overline{O}_V$ admits a $1$-step desingularization $Z$ then ${\mathbf F}_{\bullet}$ is a minimal free resolution of the normalization of $\overline{O}_V$.
\ec
\bpr If $Q$ is extended Dynkin, then $E_Q \geq 0$. This implies ${\mathbf F}_i=0$ for $i < 0$. The result then follows from Theorem \ref{thm3.3}. 
\qed

\section{Examples}\label{sec4}

\be{ex1} Consider $Q=A_4$ with orientation as in the figure below.
  \begin{figure}[h]
  \includegraphics[scale=0.5,clip]{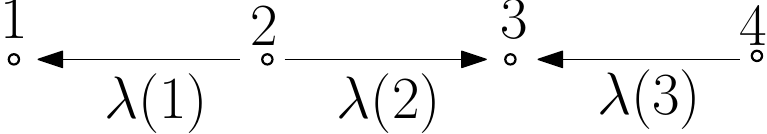}
  \label{a4}
  \caption{$A_4$}
  \end{figure}
Let $V$ be the direct sum of indecomposables with dimension vectors $(1, 0, 0, 0)$, $(1, 1, 1, 0)$, $(0, 0, 1, 0)$, $(0, 0, 1, 1)$, $(0, 1, 1, 0)$, $(1, 1, 1, 1)$ and $(1, 1, 0, 0)$. $V$ admits a $1$-step desingularization with dimension vectors $\alpha=(4, 4, 5, 2)$ and $\beta=(2, 3, 2, 1)$. The coordinate ring of $Rep(Q,(4, 4, 5, 2))$ is \[A=Sym(V_2 \otimes V_1^*)\oplus Sym (V_2 \otimes V_3^*)\oplus Sym (V_4\otimes V_3^*)\] Let $R_i$ denote the subspace of $V_i$ of dimension $\beta_i$ and let $Q_i:=V_i/R_i$. Then \[\xi=R_2 \otimes Q_1^* \oplus R_2 \otimes Q_3^* \oplus R_4 \otimes Q_3^*\]

\[\bigwedge^t \xi = \bigoplus_{\sum_{i=1}^{3}|\lambda(i)|=t} S_{\lambda(1)'} Q_1^* \otimes S_{\lambda(12)}R_2\otimes S_{\lambda(23)'}Q_3^*\otimes S_{\lambda(3)}R_4 \]

The resolution of $\overline{O}_V$ is-
\begin{center}
	     $A$  \\
	     
	     \begin{center} $ \uparrow  $ \end{center}
	     
	     $\displaystyle{(\wedge^4V_2\otimes\wedge^4V_3^*\otimes A(-7))\oplus (\wedge^4V_1^* \otimes \wedge^4V_2 \otimes A(-7))} $
	     $\displaystyle(\wedge^3V_2\otimes \wedge^5V_3^*\otimes \wedge^2V_4\otimes A(-9) \oplus(\wedge^4V_1^*\otimes S_{2221}V_2\otimes 
	     \wedge^5V_3^*\otimes\wedge^2V_4\otimes A(-17))$\\
	     
  	   \begin{center} $ \uparrow  $ \end{center}
	   	   
	   	 $\displaystyle{(S_{2111}V_2 \otimes \wedge^5V_3^* \otimes A(-8))
	   	 \oplus(\wedge^4V_1^*\otimes S_{2222}V_2 \otimes \wedge^4V_3^*\otimes A(-14))}$\\
	     $\displaystyle{(\wedge^4V_1^*\otimes S_{3222}V_2\otimes \wedge^5V_3^* \otimes A(-15))}$\\ 	
	     
	   	   \begin{center}$  \uparrow  $\end{center}  
	   	   
	     $\displaystyle{(\wedge^4V_1^*\otimes S_{3222}V_2\otimes \wedge^5V_3^* \otimes A(-15))}$\\ 	    
	     $\displaystyle{\oplus(\wedge^4V_1^* \otimes S_{2222}V_2 \otimes S_{21111}V_3^* \otimes\wedge^2V_4\otimes A(-17))}$\\ 
	     $\displaystyle{\oplus( S_{2222}V_2 \otimes S_{22222}V_3^* \otimes \wedge^2V_4 \otimes A(-17))}$\\  
	     
	     \begin{center}$  \uparrow  $\end{center}  
	   	   
	     $\displaystyle{\wedge^4V_1^*\otimes S_{3333}V_2^*\otimes S_{22222}V_3 \otimes \wedge^2V_4 \otimes A(-24)}$\\ 	    
	 \end{center} 
\ee

\vskip1cm
 
\be{ex2} Let $Q=D_5$ with the following orientation
\begin{figure}[h]
  \includegraphics[scale=0.5,clip]{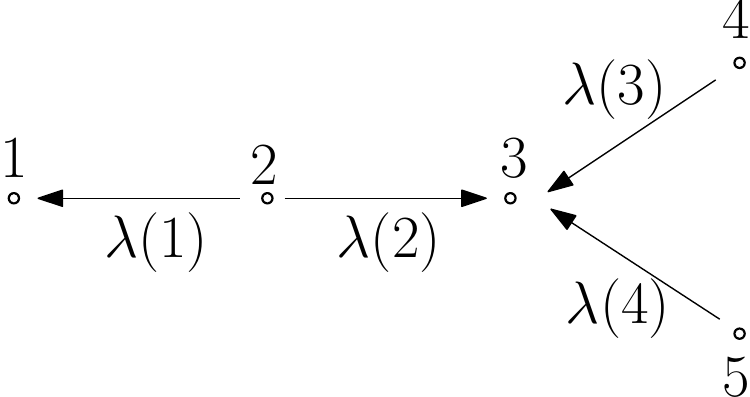}
    \label{d5}
  \caption{$D_5$}
  \end{figure}
  
Let $V$ be the direct sum of indecomposables with dimension vectors $(1, 0, 0, 0, 0)$, $(1, 1, 1, 0, 0)$, $(0, 0, 1, 0, 1)$, $(0, 0, 1, 1, 1)$ and $(1, 2, 2, 1, 1)$. $V$ admits a $1$-step desingularization with dimension vectors $\alpha=(3, 3, 5, 2, 3)$ and $\beta=(1, 2, 3, 2, 2)$. Then 
\[A=Sym(V_2 \otimes V_1^*)\oplus Sym (V_2 \otimes V_3^*)\oplus Sym (V_4 \otimes V_3^*)\oplus Sym (V_5 \otimes V_3^*)\]
 Let $R_i$ denote the subspace of $V_i$ of dimension $\beta_i$ and let $Q_i:=V_i/R_i$. Then \[\xi=R_2 \otimes Q_1^* \oplus R_2 \otimes Q_3^* \oplus R_4 \otimes Q_3^* \oplus R_5 \otimes Q_3^*\]

\[\bigwedge^t \xi= \bigoplus_{\sum_{i=1}^{3}|\lambda(i)|=t} S_{\lambda(1)'} Q_1^* \otimes S_{\lambda(12)}R_2 \otimes S_{\lambda(234)'}Q_3^*\otimes S_{\lambda(3)}R_4\otimes S_{\lambda(4)}R_5 \]

The resolution of $\overline{O}_V$ is-

\begin{center}
	     $A$\\
	     
	     \begin{center} $ \uparrow  $ \end{center}
	     
	     $\displaystyle{(\wedge^3V_1^*\otimes\wedge^3V_2 \otimes A(-5))\oplus(\wedge^3V_2\otimes\wedge^5V_3^*\otimes \wedge2V_4
	     \otimes A(-9)}$\\
	     $\displaystyle{\oplus(\wedge^5V_3^* \otimes \wedge^2V_4 \otimes \wedge^3V_5 \otimes A(-9))} $\\
	     $\displaystyle{\oplus(\wedge^2V_1^*\otimes\wedge^3V_2\otimes \wedge^5V_3^*\otimes\wedge^1V_4\otimes \wedge^3V_5\otimes A(-13))}$\\
	     $\displaystyle{\oplus(\wedge^2V_1^*\otimes\wedge^3V_2\otimes \wedge^5V_3^*\otimes\wedge^2V_4\otimes \wedge^2V_5\otimes A(-13))}$\\

	   	   \begin{center} $ \uparrow  $ \end{center}

	     $\displaystyle{\oplus(\wedge^3V_1^*\otimes\wedge^3V_2\otimes \wedge^5V_3^*\otimes\wedge^2V_4\otimes \wedge^3V_5\otimes A(-14))}$\\ 	     $\displaystyle{\oplus(\wedge^3V_1^*\otimes S_{211}V_2\otimes \wedge^5V_3^*\otimes\wedge^1V_4\otimes \wedge^3V_5\otimes A(-14))}$\\
	     $\displaystyle{\oplus(\wedge^3V_1^*\otimes S_{211}V_2\otimes \wedge^5V_3^*\otimes\wedge^2V_4\otimes \wedge^2V_5\otimes A(-14))}$\\
	     $\displaystyle{\oplus(\wedge^3V_1^*\otimes S_{222}V_2\otimes \wedge^5V_3^*\otimes\wedge^2V_4 \otimes A(-14))}$\\
	     $\displaystyle{\oplus(\wedge^2V_1^*\otimes \wedge^3V_2\otimes S_{21111}V_3^*\otimes\wedge^2V_4\otimes \wedge^3V_5\otimes 
	     A(-14))}$\\
	  	 $\displaystyle{\oplus(\wedge^3V_2\otimes S_{22222}V_3^*\otimes S_{22}V_4\otimes \wedge^3V_5\otimes A(-18))}$\\
	     $\displaystyle{\oplus(\wedge^2V_1^*\otimes S_{222}V_2\otimes S_{22222}V_3^*\otimes S_{21}V_4\otimes \wedge^3V_5\otimes A(-22))}$\\
	     $\displaystyle{\oplus(\wedge^2V_1^*\otimes S_{222}V_2\otimes S_{22222}V_3^*\otimes S_{22}V_4\otimes \wedge^2V_5\otimes A(-22))}$\\    	   	 
	   	   \begin{center}$  \uparrow  $\end{center}  
	   	  
	   	 $\displaystyle{\oplus(\wedge^3V_1^*\otimes S_{211}V_2\otimes S_{21111}V_3^*\otimes \wedge^2V_4\otimes\wedge^3V_5\otimes A(-15)}$\\ 
	   	 $\displaystyle{\oplus(\wedge^2V_1^*\otimes S_{222}V_2\otimes S_{32222}V_3^*\otimes S_{22}V_4\otimes \wedge^3V_5\otimes A(-23))}$\\
	     $\displaystyle{\oplus(\wedge^3V_1^*\otimes S_{322}V_2\otimes S_{22222}V_3^*\otimes S_{22}V_4\otimes \wedge^2V_5\otimes A(-23))}$\\
 	   	 $\displaystyle{\oplus(\wedge^3V_1^*\otimes S_{322}V_2\otimes S_{22222}V_3^*\otimes S_{21}V_4\otimes \wedge^3V_5\otimes A(-23))}$\\
	     $\displaystyle{\oplus(\wedge^3V_1^*\otimes S_{222}V_2\otimes S_{22222}V_3^*\otimes S_{22}V_4\otimes \wedge^3V_5\otimes A(-23))}$\\
     
	   	   \begin{center}$  \uparrow  $\end{center}
	   	   
	     $\displaystyle{(\wedge^3V_1^*\otimes S_{322}V_2\otimes S_{32222}V_3^*\otimes S_{22}V_4\otimes \wedge^2V_5\otimes A(-24))}$\\

	 \end{center} 
\ee

\vskip1cm
	 
\be{ex3} $Q=E_6$ with the orientation 
  \begin{figure}[h]
  \includegraphics[scale=0.5,clip]{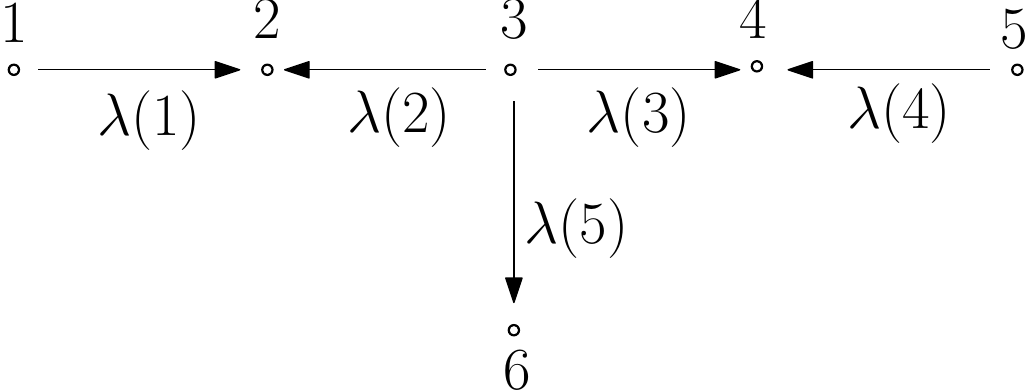}
  \label{e6}
  \caption{$E_6$}
  \end{figure}
  
Let $V=I_1 \oplus I_2$ where $I_1$ and $I_2$ are the indecomposable representations with dimension vectors $(0, 1, 1, 1, 0, 1)$ and $(1, 2, 3, 2, 1, 1)$ respectively.
Then $V$ admits a $1$-step desingularization with dimension vectors $\alpha=\underline{\hbox{dim}}~V=(1, 3, 4, 3, 1, 2)$ and $\beta=(1, 2, 3, 2, 1, 1)$. 
\[A=Sym(V_1 \otimes V_2^*)\oplus Sym (V_3 \otimes V_2^*)\oplus Sym (V_3 \otimes V_4^*)\oplus Sym (V_3 \otimes V_6^*)\oplus Sym (V_5 \otimes V_4^*) \]

\[\xi= R_1 \otimes Q_2^* \oplus R_3\otimes Q_2^* \oplus R_3 \otimes Q_4^* \oplus R_3\otimes Q_6^* \oplus R_5\otimes Q_4^*\]

\[\bigwedge^t \xi = \bigoplus_{\sum_{i=1}^{5}|\lambda(i)|=t} S_{\lambda(1)} R_1\otimes S_{\lambda(12)'}Q_2^*\otimes S_{\lambda(235)}R_3  \otimes S_{\lambda(34)'}Q_4^* \otimes S_{\lambda(4)}R_5\otimes S_{\lambda(5)'}Q_6^*\]

The resolution of $\overline{O}_V$ is
\begin{center}
	     $A$\\
	     	     
	     \begin{center} $ \uparrow  $ \end{center}
	     
	     $\displaystyle{(V_1\otimes\wedge^3V_2^* \otimes \wedge^4V_3\otimes\wedge^2V_6^*\otimes 
	     A(-9))\oplus(\wedge^4V_3\otimes\wedge^3V_4^*\otimes V_5\otimes \wedge^2V_6^*\otimes A(-9)}$\\
	     $\displaystyle{\oplus(V_1\otimes\wedge^3V_2^* \otimes \wedge^4V_3 \otimes \wedge^3V_4^* \otimes V_5 \otimes A(-11))} $
	     $\displaystyle{\oplus(V_1\otimes\wedge^3V_2^* \otimes S_{2221}V_3\otimes\wedge^3V_4^*\otimes \wedge^2V_6^*\otimes A(-15))}$\\
	     $\displaystyle{\oplus(V_1 \otimes\wedge^3V_2^*\otimes S_{222}V_3\otimes \wedge^3 V_4^*\otimes V-5\wedge^2V_6^*\otimes A(-15))}$\\ 
	     $\displaystyle{\oplus(\wedge^3V_2^* \otimes S_{2222}V_3 \otimes \wedge^3V_4^* \otimes \wedge^2V_6^* \otimes A(-15))} $\\
	     $\displaystyle{\oplus(\wedge^3V_2^* \otimes S_{2221}V_3 \otimes \wedge^3V_4^* \otimes V_5 \wedge^2V_6^* \otimes A(-15))} $\\  
	      
	   	   \begin{center} $ \uparrow  $ \end{center}
	   	   
	   	 $\displaystyle{(V_1\otimes\wedge^3V_2^* \otimes S_{2222}V_3 \otimes \wedge^3V_4^*\otimes S_{21}V_6^*\otimes A(-16))}$\\
	   	 $\displaystyle{\oplus(V_1\otimes\wedge^3V_2^* \otimes S_{2221}V_3 \otimes \wedge^3V_4^*\otimes V_5 \otimes S_{21}V_6^*\otimes 
	   	 A(-16))}$\\
	   	 $\displaystyle{\oplus(V_1\otimes\wedge^3V_2^* \otimes S_{2221}V_3 \otimes S_{211}V_4^*\otimes V_5 \otimes\wedge^2V_6^*\otimes 
	   	 A(-16))}$\\
	   	 $\displaystyle{\oplus(V_1\otimes S_{211}V_2^* \otimes S_{2222}V_3 \otimes\wedge^3V_4^*\otimes V_5\otimes\wedge^2V_6^*\otimes 
	   	 A(-16))}$\\ 
	     $\displaystyle{\oplus(V_1\otimes S_{211}V_2^* \otimes S_{2221}V_3 \otimes\wedge^3V_4^*\otimes V_5\otimes\wedge^2V_6^*\otimes 
	     A(-16))}$\\ 	    
	     $\displaystyle{\oplus(\wedge^3V_2^* \otimes S_{2222}V_3 \otimes\wedge^3V_4^*\otimes V_5\otimes S_{21}V_6^*\otimes A(-16))}$\\ 
	     $\displaystyle{\oplus(\wedge^3V_2^* \otimes S_{2222}V_3 \otimes S_{211}V_4^*\otimes V_5\otimes\wedge^2V_6^*\otimes A(-16))}$\\  
	   	 
	   	   \begin{center}$  \uparrow  $\end{center}  
	   	   
	     $\displaystyle{(V_1\otimes \wedge^3V_2^* \otimes S_{2222}V_3 \otimes S_{211}V_4^*\otimes V_5\otimes S_{21}V_6^*\otimes 
	     A(-17))}$\\ 	    
	     $\displaystyle{\oplus(V_1 \otimes S_{211}V_2^* \otimes S_{2222}V_3 \otimes\wedge^3V_4^*\otimes V_5\otimes S_{21}V_6^*\otimes 
	     A(-17))}$\\ 
	     $\displaystyle{\oplus(V_1\otimes S_{211}V_2^* \otimes S_{2222}V_3 \otimes S_{211}V_4^*\otimes V_5\otimes\wedge^2V_6^*\otimes 
	     A(-17))}$\\  
	 \end{center} 
\ee

\bibliographystyle{alpha}
\bibliography{references}

\end{document}